\documentclass[11pt]{article}
\usepackage[mathscr]{eucal}   
\usepackage{mathrsfs}
\usepackage{amsfonts}
\usepackage{amsmath}
\usepackage{amsthm}
\usepackage{amssymb}
\usepackage{latexsym}

\usepackage{pstricks,pst-node}

\theoremstyle{plain}
\newtheorem{definition}[equation]{Definition}
\newtheorem{corollary}[equation]{Corollary}
\newtheorem{lemma}[equation]{Lemma}
\newtheorem{proposition}[equation]{Proposition}
\newtheorem{theorem}[equation]{Theorem}
                                                                                
\theoremstyle{definition}
\newtheorem{remark}[equation]{Remark}
\newtheorem{example}[equation]{Example}
                                                                                
\numberwithin{equation}{section}
                                                                                
\makeatletter
\renewcommand{\subsection}{\@startsection{subsection}{2}{0pt}{-3ex
plus -1ex minus -0.2ex}{-2mm plus -0pt minus
-2pt}{\normalfont\bfseries}} \makeatother

\def\C{{\mathbb{C}}}
\def\g{{\mathfrak{g}}}

\newcommand{\Id}{{\mathtt{Id}}}
\renewcommand{\H}{{\mathscr{H}}}

\newcommand{\cc}{{\mathsf{C}}}
\newcommand{\hh}{{\mathsf{H}}}
\renewcommand{\O}{{\mathcal{O}}}
\newcommand{\Ug}{{\mathcal{U}\mathfrak{g}}}
\newcommand{\Tr}{{{\mathsf {Tr}}}}
\newcommand{\ot}{{\otimes}}
\newcommand{\ov}{{\overline{v}}}
\newcommand{\oq}{{\overline{q}}}
\renewcommand{\t}{\tau}
\newcommand{\h}{{\mathfrak{h}}}
\newcommand{\sset}{{\subset}}
\def\hp{\hphantom{x}}
\begin{document}

\centerline{\Large {\textbf{Continuous Hecke
algebras}}}
\bigskip
\centerline{\large {\sc {Pavel Etingof, Wee Liang Gan, Victor Ginzburg}}}
\bigskip

\bigskip

\centerline{\sf Table of Contents}
\vskip -1mm

$\hspace{30mm}$ {\footnotesize \parbox[t]{115mm}{
\hp${}_{}$\!\hp1.{ $\;\,$} {\tt Introduction}\newline
\hp2.{ $\;\,$} {\tt Continuous Hecke algebras}\newline
\hp3.{ $\;\,$} {\tt Continuous symplectic reflection 
algebras and Cherednik algebras}\newline
\hp4.{ $\;\,$} {\tt Infinitesimal Hecke algebras}\newline
\hp5.{ $\;\,$} {\tt Representation theory of 
continuous Cherednik algebras}\newline
\hp6.{ $\;\,$} {\tt Case of  wreath-products}
}}

\vskip .5in

{\sl The theory of PBW properties of quadratic algebras, 
to which this paper aims to be a modest contribution, 
originates from the pioneering work of Drinfeld (see \cite{Dr1}). 
In particular, as we learned after publication of \cite{EG} (to
the embarrassment of two of us!), symplectic reflection algebras, as well as PBW theorems
for them, were discovered by Drinfeld in the classical paper \cite{Dr2} 15 years
before \cite{EG} (namely, they are a special case of degenerate affine Hecke
algebras for a finite group $G$ introduced in \cite{Dr2}, Section 4). 

It is our great pleasure to dedicate this paper
to Vladimir Drinfeld on the occasion of his 50-th
Birthday.}

\vskip .5in

\section{Introduction}

In \cite{Dr2}, Section 4, Drinfeld introduced 
the family of degenerate affine (=graded) Hecke algebras
attached to any finite group $G$ and its linear finite
dimensional complex representation $V$. The goal of this paper is to
define and study the ``continuous'' generalization of 
these algebras, in which the group $G$ is a reductive
algebraic group, and $V$ is its algebraic representation. 
We call this generalization {\bf continuous Hecke algebras.}  
They include {\bf continuous symplectic reflection algebras} and 
{\bf continuous Cherednik algebras}, which are generalizations
of symplectic reflection algebras from \cite{EG} and rational Cherednik 
algebras. 

A motivation for studying continuous Hecke algebras 
comes from the fact that their representation theory
(which is yet to be developed) unifies the representation 
theories of real reductive groups, 
Drinfeld-Lusztig degenerate affine Hecke algebras, and 
symplectic reflection algebras of \cite{EG} (in particular, 
rational Cherednik algebras).  

Let us briefly describe the contents of the paper. 
Let $G$ be a reductive algebraic group
\footnote{Unless otherwise
specified, all algebraic groups and Lie algebras in this paper are over $\C$.}
with an algebraic representation $\rho: G\to GL(V)$. 
Let $\O(G)^*$ be the algebra, with respect to convolution, of algebraic
distributions on $G$ (see subsection \ref{ad}). 
This algebra has a natural $\O(G)$-module structure.

Let $\kappa\in (\O(G)^*\otimes \wedge^2
V^*)^G$. Following the idea of \cite{Dr2}, we define the algebra $\hh_\kappa$ to be the quotient 
of $TV\rtimes \O(G)^*$ by the ``main commutation relation''
$$
[x,y]=\kappa(x,y), \quad x,y\in V, 
$$
and we are interested in finding those $\kappa$ 
for which $\hh_\kappa$ is a flat deformation of $\hh_0=SV\rtimes
\O(G)^*$ (the PBW property). As in \cite{EG}, it turns out that flatness is
equivalent to the Jacobi identity
$$
\kappa(x,y)(z-z^g)+\kappa(y,z)(x-x^g)+\kappa(z,x)(y-y^g)=0,\
x,y,z\in V,\ g\in G.
$$ 
We call the algebras $\hh_\kappa$ satisfying the PBW property 
(equivalently, the Jacobi identity) {\it continuous Hecke
algebras}.

One of the main problems treated in this paper is 
to classify continuous Hecke algebras,
i.e. to find explicitly all distributions $\kappa$ satisfying the
Jacobi identity. 
In the case when $G$ is finite, the answer is given in \cite{Dr2}
(see also \cite{RS}). However, for an infinite $G$ the situation is
more complicated and we don't have a complete answer; so we
present some partial results.  

Let $\Phi\subset G$ be the closed subscheme of $G$ defined by the
equation $\wedge^3(1-g|_V)=0$. Our first result is that 
if $\kappa$ satisfies the Jacobi identity then 
it is supported on $\Phi$. 

Next, let $\theta\in(\O(\Phi)^*\otimes
\wedge^2 V^*)^G$,
and $\tau\in (\O({\rm Ker}\rho)^*\otimes \wedge^2V^*)^G$. 
We show that the distribution 
$$
\kappa(x,y):=\tau(x,y)+\theta((1-g)x,(1-g)y)
$$ 
satisfies the Jacobi identity and hence defines a continuous
Hecke algebra. This is the general answer if the group $G$ is
finite. 

Another collection of examples comes from symmetric
pairs $(\bold G,\bold K)$ of reductive groups. 
In this case $\hh_\kappa=A(\bold G,\bold
K)$, the algebra of distributions on $\bold G$ set-theoretically supported 
on $\bold K$. 

Further, we consider the special case when 
the space $V$ carries a $G$-invariant symplectic 
form $\omega$. In this case, let $\Sigma$
be the closed subscheme of $G$ defined by the equation \linebreak
$p\circ \wedge^3(1-g|_V)=0$, where 
$p:\wedge^3V\to V$ is the projection defined by $\omega$.  
Let 
$$
\kappa(x,y)=\omega(x,y)t+\omega((1-g)x,(1-g)y)c,
$$
where $t$ is a $G$-invariant distribution 
on ${\rm Ker}\rho$, and $c$ is a $G$-invariant 
distribution on $\Sigma$. 
Here $\omega((1-g)x,(1-g)y)c$ demotes the product of 
the regular function $\omega((1-g)x,(1-g)y)$ by 
the algebraic distribution $c$, using the $\O(G)$-module structure 
on $\O(G)^*$.  

 We show that $\kappa$ 
satisfies Jacobi identity, and call the corresponding 
algebras $\hh_\kappa$ {\bf continuous symplectic 
reflection algebras}. If $G$ is finite, 
and $\rho$ is a faithful representation, they coincide 
with symplectic reflection algebras defined in \cite{EG}.  

Now let us assume in addition that $V=\h\oplus \h^*$, 
where $\h$ is a representation of $G$ (and the form $\omega$ is 
given by the natural pairing between $\h$ and $\h^*$). 
In this case we show that $\kappa$ satisfies Jacobi
identity if $\kappa(y,x)=0$ if $y,x\in \h$ 
or $y,x\in \h^*$, and 
$$
\kappa(y,x)=(y,x)t+(y,(1-g)x)c,\ y\in \h, x\in \h^*. 
$$
where $t$ is a $G$-invariant distribution on ${\rm Ker}\rho$, 
and $c$ is an invariant distribution 
on the closed subscheme $\Psi\sset G$,  ``of complex reflections'', 
defined by the equation $\wedge^2(1-g|_\h)=0$.
We call algebras $\hh_\kappa$ for such $\kappa$ 
{\bf continuous Cherednik algebras}.
 If $G$ is finite and $\rho$ is a faithful representation, they
coincide with the rational Cherednik algebras considered in
\cite{EG}. 

One of our main results is that if $\h$ is irreducible 
of real or complex type and faithful then the classes of 
continuous Hecke algebras, continuous symplectic reflection
algebras, and continuous Cherednik algebras are the same.  
 
We consider in detail the examples 
$G=GL_n,O_n,Sp_{2n}$, and $V=\Bbb C^{2n}$ with the natural
symplectic form. 
\footnote{The example $G=Sp_2=SL_2$ was considered earlier in
\cite{Kh}.}
 In all these examples, any continuous Hecke 
algebra $\hh_\kappa$ is a continuous symplectic reflection
algebra, and we classify explicitly the possible distributions
$c$. In particular, for $G=O_n$ we find that there is just one 
essential parameter, while for $G=GL_n$ and $G=Sp_{2n}$ there are
infinitely many. 
 
If the distribution $\kappa$ is set-theoretically supported at
$1$ (i.e. belongs to $\Ug$, where $\g$ is the Lie algebra 
of $G$), then one can define the algebra $\H_\kappa$ to be the
quotient of $TV\rtimes \Ug$ by the main commutation relation. 
It is a flat deformation of $SV\rtimes \Ug$ 
if and only if $\kappa$ satisfies the Jacobi identity; 
in this case we call $\hh_\kappa$ an {\bf infinitesimal Hecke
algebra}. We write such $\kappa$ explicitly in the cases $G=GL_n, G=Sp_{2n}$
(in the case $G=O_n$, there are no interesting algebras
$\H_\kappa$). 

An interesting problem 
is to develop representation theory of the algebras $\hh_\kappa$ and
$\H_\kappa$. In the case of $G=SL_2$, $V=\Bbb C^2$, this was  
done in \cite{Kh}, but other examples 
(apart from finite $G$) remain mostly unexplored.
We begin to tackle this problem by
introducing the category $\O$ and studying the Dunkl-Cherednik  
representation of continuous Cherednik algebras.  

\begin{remark} Note that continuous representations of the algebra
$A(\bold G,\bold K)$ are nothing but Harish-Chandra modules 
for the pair $(\bold G,\bold K)$. Thus, the representation theory
of continuous Hecke algebras is a generalization of the 
representation theory of real reductive groups. 
\end{remark}

At the end of the paper, we consider the case 
of wreath products $G=S_n\ltimes \Gamma^n$, where $\Gamma$ is a
reductive subgroup of $SL_2$, and $V=\Bbb C^{2n}$. 
We consider in detail the three cases 
in which $\Gamma$ is infinite, and generalize the results of
Montarani on the classification of 
representations of $\hh_\kappa$ irreducible over $G$ 
(proved in \cite{M} in the case when $\Gamma$ is finite). 

{\bf Acknowledgments.}
The work of P.E. was partially supported by the NSF grant DMS-9988796.
The work of V.G. was partially supported by the NSF grant
DMS-0303465. The work of W.L.G. was partially supported by the
NSF grant  DMS-0401509. 
The work of P.E. and V.G. was partially supported by 
the CRDF grant RM1-2545-MO-03. 

\section{Continuous Hecke algebras}

\subsection{Algebraic distributions.}\label{ad}
Let $X$ be an affine scheme of finite type 
over $\C$. We shall denote by $\O(X)$ the algebra
of regular functions on $X$. An \emph{algebraic distribution} on $X$ is an
element in the dual space $\O(X)^*$ of $\O(X)$. 
For $c\in O(X)^*$ and $f\in \O(X)$, we will denote the value $c(f)$
by $(c,f)$. 

The space $\O(X)^*$ is naturally equipped 
with the weak (inverse limit) topology. 
Note also that $\O(X)^*$ is a module over $\O(X)$: for any 
$f\in \O(X)$ and $\mu\in \O(X)^*$ we can define 
the element $f\mu=\mu f\in \O(X)^*$ by $(f\mu,g)=(\mu,fg)$.

Let $Z$ be a closed subscheme of $X$, and write $I(Z)$ for its defining
ideal in $\O(X)$. We say that an algebraic distribution $\mu$ on $X$
is \emph{supported} on the scheme $Z$ if $\mu$ annihilates $I(Z)$.
Clearly, the space of algebraic distributions on $X$ supported on $Z$ is 
naturally isomorphic to the space of algebraic distributions on
$Z$. 

Now assume that $Z$ is reduced. We say that $\mu\in \O(X)^*$
is \emph{scheme-theoretically} (respectively, \emph{set-theoretically})
\emph{supported} on the set $Z$ if $\mu$ annihilates $I(Z)$ 
(respectively, some power of $I(Z)$).

\begin{example}
For each point $a\in X$, the delta function $\delta_a\in \O(X)^*$ is
defined by $\delta_a(f):=f(a)$ where $f\in \O(X)$. It
is scheme-theoretically supported at the point $a$, and its derivatives
are set-theoretically supported at this point.
\end{example}

Let $G$ be a reductive algebraic group.
Since $\O(G)$ is a coalgebra, its dual space $\O(G)^*$ is an algebra under
convolution. The unit of this algebra is the delta function
$\delta_1$ of the identity element $1\in G$. 

Note that a continuous representation of the algebra $\O(G)^*$ 
is the same thing as a locally finite $G$-module (i.e. a
$G$-module which is a direct sum of finite dimensional algebraic
representations of $G$).  

Suppose that $G$ acts on $X$.
Then $G$ acts also on $\O(X)$ and $\O(X)^*$. 
We have $\O(X) = \bigoplus_V M_V\ot V$, where $V$ runs over the
irreducible representations of $G$, and $M_V$ are multiplicity spaces. Thus,
$\O(X)^* = \prod_V M^*_V \ot V^*$.
In particular, $\O(G)^* = \prod_V V\ot V^*$ as a
$G\times G$-module.

Recall that the categorical quotient $X/G$ is an affine scheme of
finite type with
${\mathcal O}(X/G)={\mathcal O}(X)^G$. 
It is clear that the space $(\O(X)^*)^G$ of invariant 
algebraic distributions on $X$ is naturally isomorphic to the space $\O(X/G)^*$
of algebraic distributions on the quotient $X/G$. 
We will denote this space by $\cc(X)$. 
Clearly, if $Z\subset X$ is a $G$-invariant closed subscheme, then
$\cc(Z)$ is naturally a subspace of $\cc(X)$. 

\subsection{Definition and properties of continuous Hecke algebras.} 

Let $G$ be a reductive algebraic group (not necessarily
connected), and $\rho: G\to GL(V)$ 
be a finite dimensional algebraic representation
of $G$. The semi-direct product $TV\rtimes \O(G)^*$ is the algebra
generated by $x\in V$ and $\mu\in \O(G)^*$ with the relations
\begin{equation} \label{semieq}
\mu \cdot x=\sum_i v_i\cdot (v_i^*,gx)\mu
\quad \mbox{ for all } x\in V,\,\mu\in \O(G)^*,
\end{equation}
where $v_i$ is a basis of $V$ and $v_i^*$ the dual basis of $V^*$.
(Here $(v_i^*,gx)\mu$ denotes the product of the regular function 
$(v_i^*,gx)$ by the distribution $\mu$).
Given a skew-symmetric $G$-equivariant $\C$-bilinear pairing
$\kappa:V \times V \to \O(G)^*$ 
(where $G$ acts on $\O(G)^*$ by conjugation), let $\hh_\kappa=\hh_\kappa(G)$
be the quotient of $TV\rtimes \O(G)^*$ by the relations
\begin{equation} \label{defeq}
[x,y] = \kappa(x,y) \quad \mbox{ for all } x,y\in V.
\end{equation}
Define a filtration on $\hh_\kappa$ by setting
$\deg(x)=1$, $\deg(\mu)=0$ for $x\in V$, $\mu\in \O(G)^*$.
Note that $\hh_0 = SV \rtimes \O(G)^*$. There is a natural
map $\hh_0 \to gr(\hh_\kappa)$ which is a surjective graded
algebra homomorphism. We say that the \emph{PBW property} holds for
$\hh_\kappa$ if this morphism is an isomorphism.

For $x\in V$, let $x^g$ denote the result of the action of 
$g$ on $x$. 

\begin{theorem}\label{jac} The algebra $\hh_\kappa$ has 
the PBW property if and only if $\kappa$ satisfies 
the Jacobi identity
\begin{equation} \label{jacobi}
(z-z^g)\kappa(x,y) + (y-y^g)\kappa(z,x)+(x-x^g)\kappa(y,z)=0
\quad \mbox{ for all } x,y,z\in V.
\end{equation}
\end{theorem}

\begin{proof}
The algebra $\hh_0$ is Koszul. Therefore, 
by the results of \cite{BG},\cite{BGS} (see also \cite[\S 2]{EG}), 
the PBW property is equivalent to the Jacobi identity
\footnote{Strictly speaking, the result of \cite{BGS}
is stated in the situation when the degree zero part of the algebra 
is finite dimensional semisimple; however, one can check that all
arguments extend to the case when the degree zero part is
the inverse limit of finite dimensional semisimple algebras, 
in particular when the degree zero part is $\O(G)^*$.} 
 
\begin{equation} \label{jeq}
[\kappa(x,y),z] + [\kappa(z,x),y]+[\kappa(y,z),x]=0 
\quad \mbox{ for all } x,y,z\in V.
\end{equation}
By (\ref{semieq}), we have $$[z,\kappa(x,y)] = (z-z^g)\kappa(x,y) \in
\big(V\ot\O(G)\big)\ot_{\O(G)}\O(G)^*=V\otimes \O(G)^*.$$
Hence, the Jacobi identity (\ref{jeq}) is equivalent to
equation (\ref{jacobi}), as needed.
\end{proof}

\begin{definition}
If $\hh_\kappa$ has the PBW property (equivalently,
$\kappa$ satisfies the Jacobi identity) then $\hh_\kappa$ is
called a continuous Hecke algebra.  
\end{definition}

Let $\Phi\subset G$ be the closed subscheme defined by the
equation $\wedge^3(1-g|_V)=0$. Then the set of closed points of
$\Phi$ is the set $S$ of elements of $G$ such that the rank of
the operator $1-g: V\to V$ is $\le 2$.

\begin{proposition} \label{suppo}
If the PBW property holds for $\hh_\kappa$, then $\kappa(x,y)$ is supported
on the scheme $\Phi$ for all $x,y\in V$. 
(In particular, $\kappa(x,y)$ is set-theoretically supported on $S$.)
\end{proposition}

\begin{proof}
Here and below, we will 
use the following notation: for any $v_1,\ldots,v_k,x,y\in V$, 
we write $$(v_1,\ldots,v_k|x,y) \mbox{ for }
(v_1-v_1^g)\wedge\cdots\wedge(v_k-v_k^g)\kappa(x,y)
\in\big( \wedge^k V\ot
\O(G)\big)\ot_{\O(G)}\O(G)^*=\wedge^kV\otimes \O(G)^*.$$
Thus, we can rewrite the condition (\ref{jacobi}) as
\begin{equation} \label{newj}
(z|x,y)+(x|y,z)+(y|z,x)=0 \quad \mbox{ for all }x,y,z\in V.
\end{equation}

To prove the proposition, we will need two lemmas. 

\begin{lemma} \label{zuxy}
If (\ref{newj}) holds, then $(z,u|x,y)=(x,y|z,u)$ for all $x,y,z,u\in V$.
\end{lemma}
\begin{proof}
We have 
\begin{gather*}
(u,z|x,y)=-(u,x|y,z)-(u,y|z,x),\\
(z,u|x,y)=-(z,x|y,u)-(z,y|u,x).
\end{gather*}
Subtracting the first equation from the second equation, and
using 
(\ref{newj}), we get
\begin{align*}
2(z,u|x,y)= & \big( (x,z|y,u)+(x,u|z,y)\big)+\big( (y,z|u,x)+(y,u|x,z)\big)\\
= & (x,y|z,u) + (y,x|u,z) \\
= & 2(x,y|z,u).
\end{align*}
\end{proof}

\begin{lemma} \label{zuvxy}
If (\ref{newj}) holds, then $(z,u,v|x,y)=0$ for all $x,y,z,u,v\in V$.
\end{lemma}
\begin{proof}
By Lemma \ref{zuxy}, we have
\begin{gather*}
(v,z,u|x,y)=(v,x,y|z,u),\\
(u,v,z|x,y)=(u,x,y|v,z),\\
(z,u,v|x,y)=(z,x,y|u,v).
\end{gather*}
Adding these three equations, we get
$$ 3(z,u,v|x,y) = (x,y,v|z,u)+(x,y,u|v,z)+(x,y,z|u,v) = 0.$$
\end{proof}

Now we prove Proposition \ref{suppo}.
The PBW property is equivalent to condition 
(\ref{newj}). Hence, the proposition
is immediate from Lemma \ref{zuvxy}.
\end{proof}

\subsection{Examples of continuous Hecke algebras.}

Let $\Theta:=(\O(\Phi)^*\otimes
\wedge^2 V^*)^G$, and $\Theta_0\subset \Theta$ 
be the subspace of those $\theta\in \Theta$ for which 
$\theta((1-g)x,(1-g)y)=0$ for any $x,y\in V$.  
Let $\theta\in \Theta/\Theta_0$,  
and let $\tau\in (\O({\rm Ker}\rho)^*\otimes \wedge^2V^*)^G$.
Define the distribution $\kappa$ 
by 
\begin{equation}\label{thet}
\kappa(x,y)=\tau(x,y)+\theta((1-g)x,(1-g)y).
\end{equation}

\begin{theorem}
The algebra $\hh_\kappa$ for such $\kappa$ has the PBW
property.  
\end{theorem}

\begin{proof} Since $(u-u^g)\tau=0$, the Jacobi identity 
reduces to the relation 
$$
(z-z^g)\theta(x-x^g,y-y^g)+
(y-y^g)\theta(z-z^g,x-x^g)+
(x-x^g)\theta(y-y^g,z-z^g)=0.
$$
This relation holds since
$\theta$ is alternating and is supported on $\Phi$. 
We are done.   
\end{proof}

Let us analyze in more detail the special case of finite group $G$
(this analysis is given in \cite{Dr2} and in a more expanded 
form in \cite{RS}).
In this case, if $\kappa$ satisfies the Jacobi identity, then by 
Proposition \ref{suppo}, $\kappa$ is supported on the set $S$. 
Therefore, the main commutation relation 
of the algebra $\hh_\kappa$ takes the form 
$$
[x,y]=\tau(x,y)+\sum_{g\in S,\rho(g)\ne 1}\theta_g(x,y)g,
$$
where $\tau\in (\Bbb C[{\rm Ker}\rho]\otimes \wedge^2V^*)^G$, and $\theta_g$ 
is a 2-form on $V$. It follows from the Jacobi identity that 
$\theta_g(x,y)=0$ if $x,y\in {\rm Ker}(1-g|_V)$. 
Since $\theta_g$ is $g$-invariant, 
this implies that the kernel of $\theta_g$ 
contains ${\rm Ker}(1-g|_V)$, and hence $\theta_g$ is 
uniquely determined by $g$ up to scaling. 
We also see that $\kappa(x,y)$ 
is given by formula (\ref{thet}) for an appropriate $\theta$. 

Furthermore, to secure the $G$-invariance of the right hand side of the
commutation relation, we need to require that the family of forms $\theta_g$
is $G$-invariant. If $C\subset S$ is a conjugacy
class, then a nonzero invariant family $\lbrace{\theta_g, g\in
C\rbrace}$ is unique up to scaling if exists, and it exists if
and only if the rank of $1-g|_V$ is exactly 2, and the centralizer
$Z_g$ of $g\in C$ acts trivially on 
the 1-dimensional space $\wedge^2{\rm Im}(1-g)$. 
Let us call such conjugacy class admissible, and 
denote the set of admissible conjugacy classes by $S_{\rm adm}$. 
Thus we have 
$$
[x,y]=\tau(x,y)+\sum_{g\in S_{\rm adm}}\theta_g(x,y)g.
$$

Note that since $g\in Z_g$, we find that if $g$ is admissible 
then $\det(g)=1$, so $g$ is conjugate in $GL(V)$ to 
${\rm diag}(\lambda,\lambda^{-1},1,...,1)$. 
This shows that the study of degenerate Hecke algebras $\hh_\kappa$ 
essentially reduces to the case when $G\subset SL(V)$.

\begin{remark} It is explained in \cite{Dr2,RS} that 
the class of algebras $\hh_\kappa$ with PBW property 
for finite $G$ in particular contains
degenerate affine Hecke algebras for finite Coxeter groups,
introduced by Drinfeld and Lusztig.  
\end{remark} 

\subsection{Continuous  Hecke algebras and 
symmetric pairs.} \label{sympar}

Let $(\bold G,\bold K)$ be a symmetric pair of complex reductive
groups. Let $\g,\frak k$ be the Lie
algebras of $\bold G,\bold K$. Then 
we have the standard decomposition 
$\g=\frak k\oplus {\frak p}$. Let $A(\bold G,\bold K)$ 
be the algebra of algebraic distributions 
on $\bold G$ set-theoretically supported on $\bold K$. It
is generated by $\O(\bold K)^*$ and ${\frak p}$, 
with appropriate relations. 

It is clear that 
$A(\bold G,\bold K)$ 
is isomorphic to $\hh_\kappa(\bold K)$ for 
 a particular distribution $\kappa$, 
given by $\kappa(x,y)=[x,y]\in {\frak k}\subset \O(\bold K)^*$, 
where $x,y\in {\frak p}$
(here we regard elements of ${\frak k}$ as distributions on
$\bold K$ supported in the first order neighborhood of the
identity). This means that $A(\bold G,\bold K)$ is a continuous
Hecke alegbra. 

In fact, one also has the converse result. 

\begin{theorem}\label{hss}
Assume that $\bold K\subset GL(E)$, and  
$\kappa: \wedge^2E\to {\mathfrak k}\subset \O(\bold K)^*$ 
is an invariant map, such that ${\rm Ker}\kappa=0$, and  
the algebra $\hh_\kappa(\bold K)$ satisfies
the PBW property. Then $\hh_\kappa(\bold K)$
comes from a symmetric pair, as above.
Moreover, the group $\bold G$ is semisimple.    
\end{theorem} 

The proof of the theorem is based on the following lemma. 
\footnote{The first author is 
grateful to D. Vogan for explaining to him a proof of 
this lemma.} 

Let $\bold K$ be a reductive algebraic group over $\Bbb C$ 
with Lie algebra ${\frak k}$, and $E$ 
a faithful algebraic representation of $\bold K$.
Let $\kappa: \wedge^2E\to {\frak k}$ be a homomorphism of $\bold
K$-modules, such that ${\rm Ker}\kappa=0$. Let ${\frak g}={\frak k}\oplus E$.
We have a natural linear map $[,]: \wedge^2{\frak g}\to {\frak
g}$, which extends the Lie bracket on ${\frak k}$ and the natural
${\frak k}$-actions on $E$, and satisfies the condition: 
$[,]|_{E\otimes E}=\kappa$.

\begin{lemma}\label{liea}
Suppose that $[,]$ is a Lie bracket on $\g$.
Then $\g$ is a semisimple Lie algebra. 
\end{lemma} 

\begin{proof} Let $B_\g,B_{\frak k}$ denote the Killing
forms of the corresponding Lie algebras. 
We have $B_{\frak g}(a,b)=B_{\frak k}(a,b)+\Tr_E(ab)$ for $a,b\in {\frak k}$. 
Since $E$
is faithful, for any  $a\ne 0$ in the Lie algebra ${\frak k}_{\Bbb R}$
of a compact form ${\bold K}_{\Bbb R}$ of $\bold K$, we have $B_{\frak
g}(a,a)<0$. Hence the kernel of $B_\g$ is contained in $E$.
This kernel is an ideal in $\g$, so 
it is contained in the kernel of $\kappa$. 
Hence it must be zero, so $B_\g$ is nondegenerate and $\g$ is
semisimple.  
\end{proof}

\begin{proof} (of Theorem \ref{hss}). Let 
$\g={\frak k}\oplus E$. 
Since $\hh_\kappa$ satisfies PBW property, 
we have $\g\subset \hh_\kappa(\bold K)$. 
Then $\g$ is closed
under bracket, so it is a Lie algebra 
as in Lemma \ref{liea}. 
It follows from Lemma \ref{liea} that $\g$ is semisimple.
We are done.   
\end{proof}

\section{Continuous symplectic reflection 
algebras and Cherednik algebras} 

In this section we will consider the special case when the
representation $V$ has a $G$-invariant symplectic form $\omega$. 

\subsection{Continuous symplectic reflection algebras.}

Let $\Sigma\supset \Phi$ be the closed subscheme of $G$ defined by the
equation $p\circ \wedge^3(1-g|_V)=0$, where $p:\wedge^3V\to V$ is
the projection defined by contracting the first two components
using $\omega$. 

\begin{theorem} \label{pbw}
Let $t\in (\O({\rm Ker}\rho)^*)^G$ and $c\in\cc(\Sigma)$. If
\begin{equation} \label{maineq}
\kappa(x,y) = \omega(x,y)t+\omega((1-g)x,(1-g)y)c
\quad \mbox{ for all }x,y\in V,
\end{equation}
then the PBW property holds for $\hh_\kappa$.
\end{theorem}

\begin{proof} 
We shall show that (\ref{jacobi}) holds when $\kappa$ is of the form
(\ref{maineq}). First note that $(v-v^g)t=0$. Hence, we have to show that 
$f(x,y,z;g)c = 0$, where
$$ f(x,y,z;g) := 
(z-z^g)\omega(x-x^g,y-y^g) + (x-x^g)\omega(y-y^g,z-z^g)
+ (y-y^g)\omega(z-z^g,x-x^g).$$
This is clear, since $c$ is supported on $\Sigma$. 
\end{proof} 

Let $S$ be the set of symplectic reflections in $G$,
i.e. elements $s\in G$ such that the rank of the operator $1-s$ on $V$ is $\le
2$. It is easy to see that the set of closed points of $\Sigma$ is $S\cup Q$, 
where $Q$ is the subvariety of elements $g\in G$ such that 
$(1-g)^2=0$. In particular, any semisimple element $g\in
\Sigma(\Bbb C)$ belongs to $S$, and if $\lambda\ne 1$ is an eigenvalue of $g$,
then $\omega((1-g)x,(1-g)y)=(1-\lambda)(1-\lambda^{-1})\omega_g(x,y)$,
where $\omega_g(x,y)=\omega(x',y')$, $x',y'$
being the projections of $x,y$ on the image of $1-g$
along the kernel of $1-g$. Hence, 
if $G$ is a finite group and $V$ is a faithful representation
then the algebra of Theorem \ref{pbw}
is exactly the symplectic reflection algebra of \cite{EG}. 

Thus, algebras $\hh_\kappa$ with 
$\kappa$ being of the form (\ref{maineq})
are natural continuous analogs 
of symplectic reflection algebras of \cite{EG}.
We will therefore call $\hh_\kappa$ a \emph{continuous symplectic
reflection algebra}. 

\subsection{Continuous  Cherednik algebras.} 

Consider separately the special case 
$V=\h\oplus \h^*$, $G\subset GL(\h)\subset Sp(V)$, 
where $\h$ is a complex vector space, 
and the form $\omega$ is the pairing between $\h$ and $\h^*$. 
In this case, let $\Psi$ be the 
closed subscheme of $G$ defined by the equation 
$\wedge^2(1-g|_\h)=0$. 
Obviously, $\Psi\subset \Phi$, and the set of
closed points of $\Psi$ is the set $S$ of all
complex (=symplectic) reflections in $G$, i.e., 
$s\in G$ such that the rank of $1-s|_\h$ is $\le 1$.  

\begin{theorem} \label{pbwc}
Let $t\in (\O({\rm Ker}\rho)^*)^G$ and $c\in\cc(\Psi)$. If
$\kappa(x,y)=0$ for $x,y\in \h$ or $x,y\in \h^*$, and 
one has 
\begin{equation} \label{maineqc}
\kappa(x,y) = (y,x)t+(y,(1-g)x)c
\quad \mbox{ for all }x\in \h^*,\ y\in \h
\end{equation}
then the PBW property holds for $\hh_\kappa$.
\end{theorem}

\begin{proof}
We shall show that (\ref{jacobi}) holds when $\kappa$ is of the form
(\ref{maineqc}). Similarly to the proof of Theorem \ref{pbw}, we have to show that 
$h_1(x,y,z;g)c = 0$, $h_2(x,y,u;g)c=0$, where
$$ h_1(x,y,z;g) := (z-z^g)(x,y-y^g)
- (x-x^g)(z,y-y^g),
$$
$$ 
h_2(x,y,u;g):=
(u-u^g)(x,y-y^g)-(y-y^g)(x,u-u^g),$$
where $x,z\in \h^*$, $y,u\in \h$. 
It is clear that $h(x,y,z;g)$ is antisymmetric in 
$x,z$, while $h_2(x,y,u;g)$ is antisymmetric in $y,u$. 
Hence, for any $x,y,z,u$ the functions $h_1(x,y,z;g)$ 
and $h_2(x,y,u;g)$ belong to 
the ideal of the scheme $\Psi$. 
This implies that $h_1c=0$, $h_2c=0$, as desired. 
\end{proof} 

Note that if $g\in \Psi(\Bbb C)$ is not unipotent, 
and $\lambda\ne 1$ is the nontrivial eigenvalue of $g$ in $\h^*$
then $(y,(1-g)x)=(1-\lambda)(y',x')$, 
where $x',y'$ are the projections of $x,y$ on the image of $1-g$
along the kernel of $1-g$. Hence, 
if $G$ is a finite group and $\h$ is a faithful representation
then the algebra of Theorem \ref{pbwc}
is exactly the rational Cherednik algebra of \cite{EG}. 

Thus, algebras $\hh_\kappa$ with 
$\kappa$ being of the form (\ref{maineqc})
are natural continuous analogs 
of rational Cherednik algebras.
We will therefore call $\hh_\kappa$ a \emph{continuous 
Cherednik algebra}. 

\begin{theorem}\label{conpbw}
Assume that $\h$ is an irreducible
faithful representation of $G$ of real or complex type (i.e., it
has no nonzero invariant skew-symmetric form).  
If $\hh_\kappa$ is a continuous Hecke algebra, then it is a 
continuous  Cherednik algebra, i.e. $\kappa$ has the form
given in Theorem \ref{pbwc}. 
\end{theorem}

\begin{proof}
For the proof we will need the following proposition. 

\begin{proposition} \label{xyh}
If the PBW property holds for $\hh_\kappa$, then
$\kappa(x,y)=0$ whenever $x,y$ are both in $\h$, or both in $\h^*$.
\end{proposition}
\begin{proof}
Suppose $x,y\in\h$. It follows from (\ref{newj}) that
$(z|x,y)=0$ for all $z\in\h^*$. Hence, it remains to show that
$\kappa(x,y)$ annihilates the constant functions on $G$.
This follows from the identity $(\bigwedge^2 \h^*)^G=0$.
Similarly for $x,y\in\h^*$.
\end{proof}

Now let us prove Theorem \ref{conpbw}.
Let $z,x\in \h^*,y\in \h$. 
By the Jacobi identity (\ref{newj}) and 
Proposition \ref{xyh}, we have 
$$
(z|y,x)=(x|y,z).
$$
Taking the inner product of both sides with $u\in\h$, we get 
\begin{equation}\label{ja}
(u,(1-g)z)\kappa(y,x)=(u,(1-g)x)\kappa(y,z).
\end{equation}
On the other hand, we have 
\begin{equation}\label{ja1}
(u-u^{g^{-1}})\kappa(y,z)=(y-y^{g^{-1}})\kappa(u,z);
\end{equation}
this follows from the Jacobi identity by applying $g^{-1}$. 
Let us now move $1-g$ to the left in the right hand side of 
(\ref{ja}), and apply (\ref{ja1}). 
This gives
\begin{equation}\label{ja2}
(u,(1-g)z)\kappa(y,x)=(y,(1-g)x)\kappa(u,z).
\end{equation}

Let $e_i$ be a basis of $\h$, $e_i^*$ the dual basis of $\h^*$. 
Let $t_{ij}=(e_i,(1-g)e_j^*),c_{ij}=\kappa(e_i,e_j^*)$. 
Then (\ref{ja2}) implies 
\begin{equation}\label{ja3}
t_{ij}c_{pq}-t_{pq}c_{ij}=0. 
\end{equation}

Let us regard all our distributions as distributions on $GL(\h)$. 
We claim that $c_{ij}=t_{ij}b$ for some distribution
$b$ on $GL(\h)$. Indeed, by (\ref{ja3}), 
$\lbrace{c_{ij}\rbrace}$ 
is a 1-cocycle in the Koszul complex of the $\O(GL(\h))$-module
$\O(GL(\h))^*$ for the sequence of elements $\lbrace{t_{ij}\rbrace}$. 
But this complex is exact, since $\O(GL(\h))^*$ is dual to a free
module, and the sequence $\lbrace{t_{ij}\rbrace}$ is regular. 
Thus, $\lbrace{c_{ij}\rbrace}$ is a coboundary
of some element $b$, as desired. 

Thus we have 
$$
\kappa(x,y)=(y,(1-g)x)b
$$
The distribution $b$ is defined up to
adding a multiple of $\delta_1$, so it is invariant under $G$.
 
Since $\kappa$ is supported on $G$, we find 
that $(y,(1-g)x)b$ is supported on $G$ for all $x,y$. 
Thus, for any function $\phi$ vanishing on $G$,
$b\phi=\gamma(\phi)\delta_1$, 
where $\gamma$ is a scalar-valued 
$G$-invariant linear functional such that 
$\gamma(\phi f)=\gamma(\phi)f(1)$, $f\in \O(GL(\h))$. 
This means that $\gamma$ has
the form $\gamma(\phi)=<(d\phi)(1),Y>$, where $Y\in{\frak{gl}}(\h)$
is an invariant element (i.e., $Y=t\cdot {\rm Id}$ is a scalar
operator, since $\h$ is irreducible). Let $c=b-Y$. Then we get 
$$
\kappa(x,y)=t(y,x)\delta_1+(y,(1-g)x)c,
$$
and $c$ is supported on $G$. 
 
It remains to show that $c$ is (scheme-theoretically) supported on
$\Psi$. But this is clear from the Jacobi identity:
$$
c(y,(1-g)x)(1-g)z=c(y,(1-g)z)(1-g)x. 
$$
The theorem is proved. 
\end{proof}

\begin{corollary}\label{sra=cher}
If $\h$ is a faithful irreducible representation of real or
complex type, then any continuous symplectic reflection algebra is a continuous 
Cherednik algebra, and vice versa.
\end{corollary}

\begin{proof} It follows from Theorems \ref{pbw}, \ref{conpbw} that 
any continuous symplectic reflection algebra is a continuous
 Cherednik algebra. To prove the converse, 
let $e_i$ be a basis of $\h$ and $e_i^*$ the dual basis of
$\h^*$. Then on $\Psi$ we have 
$$
\omega((1-g)x,(1-g)y)=\sum_i ((1-g)x,e_i)(e_i^*,(1-g)y)=
-\sum_i (gx,(1-g)e_i)(e_i^*,(1-g)y)=
$$
$$
-\sum_i (gx,(1-g)y)(e_i^*,(1-g)e_i)=
\sum_i ((1-g)x,y)(e_i^*,(1-g)e_i)=\phi(g)(y,(1-g)x),
$$
where $\phi(g)=\Tr(1-g)|_\h$. 
Now let $c$ be a $G$-invariant distribution on
$\Psi$, and $\widetilde {c}$ be a $G$-invariant distribution 
on $G$ such that $\phi \widetilde c=c$ (it exists 
since $\phi$ is not a zero divisor in $\O(G)^G$, as for any 
element $g\ne 1$ of a maximal compact subgroup of $G$, one has 
$\phi(g)>0$), so 
the multiplication operator by $\phi$ is injective on $\O(G)^G$ and hence 
surjective on $(\O(G)^G)^*$). 
Then we have 
$$
\omega((1-g)x,(1-g)y)\widetilde{c}=(y,(1-g)x)c.
$$
 From the Jacobi identity it is clear that $\widetilde{c}$ is
supported on $\Sigma$, which implies the required statement. 
\end{proof} 

\begin{corollary} 
Let $(\bold G,\bold K)$ be a Hermitian symmetric pair
(with semisimple $\bold G$).
Then $A(\bold G,\bold K)$ is a continuous Cherednik algebra.
\end{corollary} 

\begin{proof} 
Let $\g,{\mathfrak k}$ be the Lie algebras 
of $\bold G,\bold K$, 
Then $\g={\mathfrak k}\oplus {\mathfrak p}_+\oplus {\mathfrak
p}_-$, and we have a linear representation $\bold K\to GL({\mathfrak p}_+)$
Also, we have $[{\mathfrak p}_+,{\mathfrak p}_+]=
[{\mathfrak p}_-,{\mathfrak p}_-]=0$, while 
$[{\mathfrak p}_+,{\mathfrak p}_-]\subset {\mathfrak k}$.
It follows (after a simple computation) that 
$A(\bold G,\bold K)=\hh_\kappa(\bold K)$, where $\kappa$ is such
that $t=0$ and $c=\Delta \delta_1$, 
$\Delta$ being the Laplacian on $\bold K$ attached to the Killing
form of $\bold G$ restricted to $\bold K$. 
\end{proof} 

{\bf Remark.} We see that 
it is possible that the scheme $\Psi/G$ is
not reduced (i.e., $\Psi/G\ne
S/G$), and hence $\cc(S)\ne \cc(\Psi)$.
Indeed, for Hermitian symmetric pairs 
one often has $S=\lbrace{ 1\rbrace}$, and yet one has a  
distribution $\kappa(x,y)$ which 
satisfies the Jacobi identity and 
is not a multiple of $\delta_1$. 

\subsection{Examples.} 

\subsubsection{The case $G=GL_n$.}\label{GL}

Let $\h$ be an $n$ dimensional vector space, and
let $G=GL_n$ be the general linear group of $\h$.
In this case, 
we have an isomorphism $S/G \to \C^*: s\mapsto \lambda=\Tr|_\h(s)-n+1$.
Thus, $\cc(S)$ is the space of Fourier series $\sum_{m\in\mathbb Z}
c_m \lambda^m$.
Also, the scheme $\Psi$ 
is well known to be reduced, so we have $\Psi=S$. 
Thus in this case $\cc(\Psi)=\cc(S)$.

\subsubsection{The case $G=O_n$.}\label{O}

Let $\h$ be an $n$ dimensional vector space
with a nondegenerate symmetric bilinear form, and
let $G=O_n$ be the orthogonal group of $\h$.
Note that for $n=2$ this example can be regarded as a continuous
limit of symplectic reflection algebras of finite dihedral groups.

In this case $S$ consists of 
$1$ and the conjugacy class of orthogonal reflections, 
so $S/G$ consists of two points and hence 
the space $\cc(S)$ is two dimensional. However, 
it turns out that $\Psi$ is not reduced in this case. 
Luckily, we have the following proposition. 

\begin{proposition}\label{O_n} For $G=O_n$, $\cc(\Psi)=\cc(S)$. 
\end{proposition}

\begin{proof} First we claim that $\Psi$ is reduced 
near every orthogonal reflection. Indeed, let 
$s={\rm diag}(-1,1,...,1)\in S$. 
Then the tangent space $T_sS$ is
the space $[s,\g]$ of symmetric matrices $A=(a_{ij})$ such that 
$a_{ij}=0$ unless $i\ne j$ and $i=1$ or $j=1$.
Now assume that $s+\hbar A$ is a point of $\Psi$ 
over the ring of dual numbers $\Bbb C[\hbar]/\hbar^2$. 
The nontrivial minors of $1-s-\hbar A$
are $\Delta_{1i,1j}=-2\hbar a_{ij}$, so we get $a_{ij}=0$ if
$i,j>1$. We also have $\det(s+\hbar A)=-1+\hbar a_{11}$, so
$a_{11}=0$. Hence, $A\in T_sS$, and thus $T_s\Psi=T_sS$. 

It remains to consider the component $\Psi_1$ of the identity in
$\Psi$. Writing $g=e^A$, where $A$ is skew-symmetric, we 
may ($G$-equivariantly) identify $\Psi_1$ with the subscheme 
of ${\frak {so}}_n$ defined by the equation $\wedge^2A=0$. 

We claim that the ring $\O(\Psi_1)$ is spanned by 
functions $1$ and $a_{ij}$, $i<j$, with trivial 
multiplication. Indeed, 
it is easy to see that $1, a_{ij}$ are 
linearly independent, and $a_{ij}a_{ki}=a_{ij}a_{ki}-a_{ii}a_{kj}=0$. 
Also, $a_{ij}a_{kl}$ is symmetric under the transposition
of $j,l$ but antisymmetric under the transposition of $i,j$ and
$k,l$, so it is zero.

Thus, $\O(\Psi_1)=\Bbb C\oplus \wedge^2\h$ as a 
$G$-module, and hence $\O(\Psi_1/G)=\Bbb C$, as desired.      
\end{proof}

We conclude that for $G=O_n$ the continuous Cherednik
algebras form a 2-parameter family (out of which only one is
essential, because of scaling).

\subsubsection{The case $G=Sp_{2n}$.}\label{SP}

Let $V$ be a $2n$ dimensional vector space
with a nondegenerate skew-symmetric bilinear form, and
let $G=Sp_{2n}$ be the symplectic group of $V$.
We have an isomorphism $S/G \to \C: s\mapsto 
w=\lambda+\lambda^{-1}=\Tr|_V(s)-2n+2$.
Thus, $\cc(S)$ is the space of Fourier series $\sum_{m\in\mathbb Z}
c_m \lambda^m$ with $c_m=c_{-m}$.

\begin{theorem}\label{sp} We have:
(i) $\Phi=S$, and $S$ is an irreducible affine algebraic
variety;

(ii) if the PBW property holds for $\hh_\kappa$, 
then $\kappa(x,y)$ is of the form (\ref{maineq}), with $c\in \cc(\Phi)=\cc(S)$.
\end{theorem}

\begin{proof} 
(i) We first show that $\Phi$ is reduced (i.e., coincides with $S$). 
Let $\overline{S}\sset S$ be the subset of elements which 
are not unipotent. Let $Gr(2n-2,2n)$ be the Grassmannian variety
of $2n-2$ dimensional subspaces in $V$.
The map $\overline{S}\to Gr(2n-2,2n): g\mapsto \mathtt{Ker}(1-g)$
is a smooth submersion, with open image (the set of nondegenerate 
codimension 2 subspaces of the symplectic space $V$). 
Hence $\overline{S}$ is smooth.

Any conjugacy class in $\overline{S}$ contains a representative
of the form
$$ g= \left(\begin{array}{ccccc}
a & b & & &  \\ 
c & d & & &  \\
& & 1 & & \\
& & & \ddots & \\
& & & & 1 \end{array}\right).$$
We claim that the tangent spaces $T_g \Phi$ and $T_g S$ are equal.
Indeed, let $g+\hbar Q$ $(Q=(q_{ij}))$ be a point of $\Phi$ 
over the ring $\Bbb C[\hbar]/\hbar^2$. 
Then we have $\Tr(g^{-1}Q)=0$, and the equations $\wedge^3(g+\hbar Q-1)=0$ 
imply that $q_{ij}=0$ for all $i,j \geq 3$. 
On the other hand, $T_g S$ is spanned by matrices of the form
$[g,M]$ where $M\in \mathfrak{sp}(2n)$, and matrices $Q=gP$, 
where $P$ has the form  
$$ \left(\begin{array}{ccccc}
\alpha & \beta & &  &   \\
\gamma & \delta & & &  \\
 & & 0 & & \\
 & & & \ddots & \\
 & & & & 0 \end{array}\right) \mbox{ where }\alpha+\delta=0. $$
Therefore, $T_g \Phi = T_g S$, and so $\Phi$ is reduced at the formal 
neighborhood of each $g\in \overline{S}$.

The set of points of $S$ where $\Phi$ is not reduced is
$G$-invariant. Therefore, to prove that this set is empty, 
it remains to show that $\Phi$ is reduced at the formal neighborhood
of $1$. Let $A$ be the algebra of functions on the formal neighborhood
of $1$ in the scheme $\wedge^3(1-g)=0$, where $g\in G$. 
Let $A_0$ be the algebra of functions on the formal neighborhood of $0$ in 
the scheme $\wedge^3 M =0$, where $M\in \mathfrak{sp}(2n)$.
Explicitly, $A=\C[[1-g]]/(\wedge^3(1-g))$ and
$A_0=\C[[M]]/(\wedge^3 M)$, where $g=(g_{ij})$ and $M=(m_{ij})$
are matrices of formal variables.
The pullback of the exponential map gives an isomorphism
$\exp^*:\C[[1-g]] \to \C[[M]]$: $g \mapsto 1+M+M^2/2!+\cdots$.
Now $\exp^*:\wedge^3(g-1)\mapsto 
\wedge^3\big(M(1+M/2! + M^2/3! + \cdots)\big)$.
The element $1+M/2! + M^2/3! + \cdots$ is invertible, 
and hence $\exp^*$ gives an isomorphism from $A$ to $A_0$.
By \cite[Theorem 1]{Ku}, $A_0$ is reduced, so $A$ is also reduced.

The irreducibility of $\Phi=S$ follows from the easy fact 
from linear algebra that 
$\overline{S}$ (which is irreducible) is a dense subset of $S$. 

(ii) Let $S'$ be the subset of elements of $S$ which 
are not unipotent and do not have eigenvalue $-1$. 
By part (i), $S$ is irreducible, 
and $S'$ is an affine open set in $S$, defined by the
condition $f^2-4\ne 0$, where 
$f(s)=\Tr(s)-2n+2$. Hence $\O(S)\subset \O(S')$. 
Thus, we have a surjective $G$-invariant map $\xi: \O(S')^{*}\to
\O(S)^*$. 

The group $G$ acts on $S'$, with $S'/G=\Bbb C\setminus
\lbrace{2,-2\rbrace }$, and the stabilizer of 
each point is isomorphic to $K=Sp(2n-2)\times \Bbb C^*$. 
Thus, $S'$ is a fibration over $\Bbb C\setminus
\lbrace{2,-2\rbrace}$
whose fiber is $G/K$, and the action of $G$ is fiberwise. 

Let $\kappa$ be a distribution satisfying the PBW property. 
Then, as follows from Proposition \ref{suppo}
and part (i), $\kappa\in (\O(S)^*\otimes \wedge^2V^*)^G$. 
Let $\overline\kappa$ be a lifting 
of $\kappa$ in $(\O(S')^*\otimes \wedge^2V^*)^G$.
Since $S'$ is fibered over $\Bbb C\setminus
\lbrace{2,-2\rbrace}$ with fiber $G/K$, 
we may view $\overline{\kappa}$ 
as a distribution on 
$\Bbb C\setminus\lbrace{2,-2\rbrace}$ with values in 
the space $(\wedge^2V^*)^K$, which for $n>1$ is a two dimensional space
spanned by $\omega(x,y)$ and $\omega((1-g)x,(1-g)y)$
(where $g$ is a semisimple element in $S$ whose centralizer
in $G$ is $K$).
In other words, we
have 
$$
\overline{\kappa}(g)(x,y)=\alpha(g)\omega(x,y)+\beta(g)\omega((1-g)x,(1-g)y), 
$$
where
$\alpha,\beta$ are $G$-invariant distributions on $S'$.
The Jacobi identity for $\kappa$ implies that 
$$
(z-z^g)\overline{\kappa}(x,y)+(x-x^g)\overline{\kappa}(y,z)+
(y-y^g)\overline{\kappa}(z,x)=0
$$
modulo distributions that vanish on $\O(S)$. 
This imposes no condition on $\beta$, and 
the condition on $\alpha$ is that the distribution
$$
\alpha(g)((z-z^g)\omega(x,y)+(x-x^g)\omega(y,z)+(y-y^g)\omega(z,x))
$$
vanishes on $\O(S)$. This implies that $\alpha$ is 
scheme-theoretically supported at
$1$, so $\alpha=t\delta_1$. The result follows. 

If $n=1$, then the proof is even simpler, since in this case
$(\wedge^2V^*)^K$ is a 1-dimensional space spanned by $\omega$. 
\end{proof}

\section{Infinitesimal Hecke algebras} \label{ug}

\subsection{The main theorem.}
Let $\g$ be the Lie algebra of $G$. The universal enveloping
algebra $\Ug$ is naturally isomorphic to the subalgebra of $\O(G)^*$ 
consisting of all algebraic distributions set-theoretically
supported at the identity element $1$ of $G$, cf. \cite[II, \S6]{DG}.
When the image of $\kappa:V\times V\to \O(G)^*$ lies in $\Ug$,
let $\H_\kappa=\H_\kappa(\g)$ be the quotient of
$TV\rtimes \Ug$ by the relations (\ref{maineq}).

We define a filtration on $\H_\kappa$ by setting $\deg(x)=1$, $\deg(u)=0$ 
for $x\in V$, $u\in\Ug$.
There is a natural map $\H_0 =SV\rtimes\Ug
\to gr(\H_\kappa)$ which is a surjective 
graded algebra
homomorphism, and we say that
the PBW property holds for $\H_\kappa$ if this map is an isomorphism.
Since $\hh_\kappa(G)=\O(G)^*\otimes_{\Ug}\H_\kappa(\g)$, 
the PBW property holds for $\H_\kappa(\g)$ if and 
only if it holds for $\hh_\kappa(G)$, i.e. when $\kappa$
satisfies the Jacobi identity. 
Whenever this is the case, we call $\H_\kappa$ 
an {\it infinitesimal Hecke algebra.}

\begin{example}
For instance, we can define infiniteimal Hecke algebras 
$\H_\kappa$ in the examples 
attached to symmetric pairs, see Section \ref{sympar};
in this case it is easy to see that $\H_\kappa$ is the universal 
enveloping algebra of the Lie algebra of the group $\bold G$. 
\end{example}

The aim of this section is to prove the following result
when $\g$ is $\mathfrak{gl}_n$ or $\mathfrak{sp}_{2n}$.
\begin{theorem} \label{ugpbw}
The PBW property holds for $\H_\kappa(\g)$ 
(i.e., it is an infinitesimal Hecke algebra) if and only if 
the relations (\ref{defeq}) are of the form 
(\ref{uglneq}) when $\g=\mathfrak{gl}_n$, and 
(\ref{uspveq}) when $\g=\mathfrak{sp}_{2n}$.
\end{theorem}
The relations (\ref{uglneq}) and (\ref{uspveq}) will be given
below.

\subsubsection{The $\g=\mathfrak{gl}_n$ case.}
Let $\h:=\C^n$, $V:=\h\oplus\h^*$ and
$\g:=\mathfrak{gl}_n=\mathfrak{gl}(\h)$.
We shall identify $\g^*$ with $\g$ via the pairing
$\g \times \g \to \C: (A,B)\mapsto \Tr(AB)$.
We shall also identify $S\g$ with $\Ug$ via the symmetrization map.
For any $x\in \h^*$, $y\in \h$, $A\in\g$, we may write
\begin{equation} \label{gln}
(x, (1-\t A)^{-1}y)\det(1-\t A)^{-1} = r_0(x,y)(A)+r_1(x,y)(A)\t
+r_2(x,y)(A)\t^2+ \cdots 
\end{equation}
where $r_i(x,y)\in S\g$ is a polynomial function on $\g$.

For each polynomial $\beta=\beta_0+\beta_1\t+\beta_2\t^2+\cdots\in\C[\t]$,
define the algebra $\H_\beta=\H_\beta(\mathfrak{gl}_n(\C))$ to be
the quotient of $T(\h\oplus \h^*)\rtimes \Ug$ by the relations
\begin{equation} \label{uglneq}
[x,x']=0,\quad [y,y']=0,\quad
[y,x]= \beta_0r_0(x,y)+\beta_1r_1(x,y)+\beta_2r_2(x,y)+\cdots,
\end{equation}
for all $x,x'\in \h^*$ and $y,y'\in \h$.

\begin{remark}
To expand the left hand side of (\ref{gln}) as a power series in $\t$,
we write
\begin{gather*}
(1-\t A)^{-1} = 1+\t A+\t^2A^2+\cdots,\\
\det(1-\t A)^{-1} = e^{-\log\det(1-\t A)} = e^{-\Tr\log(1-\t A)}
= e^{\sum_{k> 0} \t^k \Tr(A^k/k)}  \\
= 1 + \t\Tr(A) + \t^2\big(\Tr(\frac{A^2}{2})
+\frac{(\Tr A)^2}{2}\big)+\cdots.
\end{gather*}
\end{remark}

\begin{example}
We have $r_0(x,y)=(x,y)$. Thus, $\H_1(\mathfrak{gl}_n(\C)) = 
D(\C^n) \rtimes \mathcal{U}\mathfrak{gl}_n(\C)$ (where $D(\Bbb
C^n)$ is the algebra of polynomial differential operators on
$\Bbb C^n$). 
\end{example}

\begin{example}
We have $r_1(x,y)=y\otimes x + (x,y)\Id\in\mathfrak{gl}_n(\C)$. 
Thus, $\H_\t(\mathfrak{gl}_n(\C)) 
= \mathcal{U}\big(\h \oplus \h^* \oplus \mathfrak{gl}_n(\C)\big)$, 
where the Lie bracket on $\h \oplus \h^* \oplus \mathfrak{gl}_n(\C)$ is
$$[(y,x,A), (y',x',A')] = (Ay'-A'y, xA'-x'A, [A,A']+r_1(x',y)-r_1(x,y')).$$
There is a Lie algebra isomorphism $\varphi: \h \oplus \h^* \oplus \mathfrak{gl}_n(\C) \to 
\mathfrak{sl}_{n+1}(\C)$ defined by 
\begin{align*}
y  \mapsto &  \left(\begin{array}{ccc|c}
                  & & &  \\
                  & 0 & & y \\
                  & & &  \\
                  \hline
                  & 0 & & 0 \end{array}\right) \quad\mbox{ if }y\in \h,\\
x  \mapsto & \left(\begin{array}{ccc|c}
                  & & &  \\
                  & 0 & & 0 \\
                  & & &  \\
                  \hline
                  & -x & & 0 \end{array}\right) \quad
\mbox{ if }x\in \h^*,\\
A \mapsto & \left(\begin{array}{ccc|c}
                  & & &  \\
                  & A & & 0 \\
                  & & &  \\
                  \hline
                  & 0 & & 0 \end{array}\right) \quad
\mbox{ if }A\in \mathfrak{sl}_n(\C) \subset\mathfrak{gl}_n(\C),\\
\lambda \Id \mapsto & \left(\begin{array}{ccc|c}
          \frac{\lambda}{n+1} & & &  \\
                  & \ddots & & 0 \\
                  & &  \frac{\lambda}{n+1} &  \\
                  \hline
                  & 0 & & \frac{-n\lambda}{n+1}\end{array}\right) \quad
\mbox{ if } \lambda\in \C.
\end{align*}

Obviously, this example comes from the hermitian symmetric pair 
$(SL_{n+1},GL_n)$, as explained in Section \ref{sympar}. 
\end{example}
                                                                                
\subsubsection{The $\g=\mathfrak{sp}_{2n}$ case.}
Let $V=\Bbb C^{2n}$ be a symplectic vector space
with symplectic form $\omega$, and $\g:=\mathfrak{sp}(V)$.
Again, we identify $\g^*$ with
$\g$ via the pairing $\g\times\g\to\C: (A,B)\mapsto\Tr(AB)$,
and $S\g$ with $\Ug$ via the symmetrization map.
For any $x,y\in V$, $A\in\g$, we may write
\begin{equation} \label{spv}
\omega(x,(1-\t^2A^2)^{-1}y)\det(1-\t A)^{-1} 
 = \ell_0(x,y)(A)+\ell_2(x,y)(A)\t^2+\ell_4(x,y)(A)\t^4 +\cdots  
\end{equation}
where $\ell_i(x,y)\in S\g$ is a polynomial function on $\g$.

For each polynomial $\beta=\beta_0+\beta_2\t^2+\beta_4t^4+\cdots\in\C[\t]$,
define the algebra $\H_\beta=\H_\beta(\mathfrak{sp}(V))$ to be
the quotient of $TV\rtimes \Ug$ by the relations
\begin{equation} \label{uspveq}
[x,y]= \beta_0\ell_0(x,y)+\beta_2\ell_2(x,y)+\beta_4\ell_4(x,y)+\cdots,
\end{equation}
for all $x,y\in V$.

\begin{remark}
Note that $A\in\mathfrak{sp}(V)$ implies $A^k=(-1)^kJ^T(A^T)^kJ$,
where $J$ is the matrix of the form $\omega$ in a symplectic basis.
Thus, in (\ref{spv}), $\omega(x, (1-\t^2A^2)^{-1}y)$ is skew-symmetric in $x,y$.
We also have $\Tr(A^k)=0$ if $k$ is odd, so
the expansion of $\det(1-\t A)^{-1}$ contains only even powers of $\t$.
\end{remark}

\begin{example}
We have $\ell_0(x,y)=\omega(x,y)$. Thus, $\H_1(\mathfrak{sp}(V))
= A_n\rtimes \mathcal{U}\mathfrak{sp}(V)$, where $A_n$ is the $n$-th Weyl 
algebra.
\end{example}

\begin{example}
When $n=1$, we have $\g=\mathfrak{sl}_2$. Let
$$A=\left(\begin{array}{rr} a & b \\ c & -a \end{array}\right),
\ e=\left(\begin{array}{rr} 0 & 1 \\ 0 & 0 \end{array}\right),
\ h=\left(\begin{array}{rr} 1 & 0 \\ 0 & -1 \end{array}\right),
\ f=\left(\begin{array}{rr} 0 & 0 \\ 1 & 0 \end{array}\right).$$
Note that $A^2=(a^2+bc)\Id$, so $\ell_{2k}(x,y)$ is equal to a constant times
$\omega(x,y) \Delta^k$, where $\Delta=\frac{h^2}{2}+ef+fe$ is the Casimir
element of $\mathcal{U}\mathfrak{sl}_2$. 
The algebra $\H_\beta(\mathfrak{sl}_2)$ was studied by Khare in \cite{Kh}.
\end{example}

\subsection{Proof of Theorem \ref{ugpbw}.}
Let $G=GL_n$ when $\g=\mathfrak{gl}_n$,
and let $G=Sp_{2n}$ when $\g=\mathfrak{sp}_{2n}$.
By Theorems \ref{pbw}, \ref{sp}, \ref{conpbw}, and Corollary \ref{sra=cher},
the PBW property holds if and only if $\kappa$ is of the form (\ref{maineq}).
Thus we have to find all $c\in\cc(S)$ such that
$\kappa(x,y)$ is set-theoretically supported at $1$ 
for all $x,y\in V$, and then compute $\kappa(x,y)$ as
an element in $\Ug$.

{\it The $\g=\mathfrak{gl}_n$ case}: 
By Proposition \ref{xyh}, we may assume that $\kappa(x,y)=0$
if $x,y$ are both in $\h$ or both in $\h^*$.
Following Subsection \ref{GL}, we identify $\cc(S)$ with the space
of algebraic distributions on $\C^*$, so that
$c\in \C[\lambda,\frac{1}{\lambda}]^*$. Then $\kappa(x,y)$ is set-theoretically 
supported at $1$ for all $x,y\in V$ if and only if
$c=c(\lambda)$ is a finite linear combination of $\delta_1$ and its
derivatives.

Now let $\h=\C^n$ be equipped with the Hermitian inner product
$\langle x,y\rangle := x^T\overline{y}$.
We may identify $\h^*$ with $\overline{\h}$ using this
inner product on $\h$. Let
\begin{gather*}
K := U(n,\Bbb C) = \{g\in G,|\, g^{-1}=\overline{g}^T\},\\
s_\theta := \mathrm{diag}(e^{i\theta}, 1, \ldots, 1) \in K,
\quad \theta\in [0,2\pi],\\
S_\theta := \{ gs_\theta g^{-1} \,|\, g\in K\} \subset K.
\end{gather*}
For any $s\in S_\theta$, denote by $p_s\in\g$ the orthogonal
projection of $\h$ to the eigenspace of $s$
with eigenvalue $e^{i\theta}$.
We have $p_s = v\ot\ov$ where
$v\in \h$ is an eigenvector of $s$ with eigenvalue $e^{i\theta}$
and $|v|=1$.
We have to compute
\begin{align*}
\kappa(x,y)= & t(x,y)+ \int_{0}^{2\pi} c(\theta) 
\Big(\int_{S_\theta} (1-e^{-i\theta})
(x,p_s y)s\,ds \Big)\,d\theta \\
= & t(x,y)+ \int_{|v|=1} (x, (v\ot\ov)y) \Big(\int_0^{2\pi}
(1-e^{-i\theta})
c(\theta)e^{i\theta (v\ot\ov)} \,d\theta\Big)\,dv \\
= & t(x,y)+ \int_{|v|=1} (x, (v\ot\ov)y) f(v\ot\ov) \,dv \ \in \Ug,
\end{align*}
where $x\in \h^*$, $y\in \h$, 
$c(\theta)$ is a finite linear combination of $\delta_0$ and its
derivatives, and $f$ is some polynomial depending on $c$.

Recall that we identify $\Ug$ with $S\g$ via the symmetrization map.
Let $F_m\in S\g$ be the function on $\g$ defined by 
$$ F_m(A) := \int_{|v|=1} \langle Av, v\rangle^{m+1}\,dv,
\quad A\in\g.$$
Note that 
$$dF_m|_A(y\ot x) = (m+1)
\int_{|v|=1} (x, (v\ot\ov)y) \langle Av, v\rangle^m\,dv.$$
Thus $\kappa(x,y)$ corresponds, under symmetrization map, 
to a finite linear combination of functions 
$dF_m|_A(y\ot x)$. 

We now compute $F_m(A)$. First, note that
\begin{align*}
\int_\h \langle Av, v\rangle^{m+1} e^{-\zeta\langle v,v\rangle}\,
dv = & \int_0^\infty e^{-\zeta r^2} \int_{|v|=r}
 \langle Av, v\rangle^{m+1} \,dv\,dr \\
= & \int_0^\infty r^{2m+2n+1} e^{-\zeta r^2}\,dr\cdot F_m(A) \\
= & \frac{(m+n)!}{2}\zeta^{-m-n-1} F_m(A).
\end{align*}
On the other hand, we have
\begin{align*}
\sum_{m=-1}^\infty \frac{1}{(m+1)!}\int_\h \langle Av, v\rangle^{m+1} 
e^{-\zeta\langle v,v\rangle}\,dv 
= & \int_\h e^{\langle Av, v\rangle - \zeta\langle v,v\rangle}\,dv\\
= & \int_\h e^{-\langle (\zeta-A)v, v\rangle}\,dv \\
= & \pi^n \det(\zeta-A)^{-1}.
\end{align*}
Hence, $F_m(A)$ is equal to a constant times the
coefficient of $\t^{m+1}$ in the expansion of $\det(1-\t A)^{-1}$
as a power series in $\t=\zeta^{-1}$.
Differentiating $\det(1-\t A)^{-1}$ and evaluating at $B\in \g$, 
we obtain
$$ \frac{\partial}{\partial B} \det(1-\tau A)^{-1}=
\frac{\Tr(\t B (1-\t A)^{-1})}{\det(1-\t A)}.$$
Setting $B=y\ot x$ gives
$\t(x,(1-\t A)^{-1}y) \det(1-\t A)^{-1}$, as desired.

\bigskip

{\it The $\g=\mathfrak{sp}_{2n}$ case}:
We let $V=\C^{2n}$ be equipped with
the symplectic form $\omega(x,y) = x^TJy$ where
$$J:= \left( \begin{array}{ccccc}
0 & 1 &  & & \\
-1 & 0 & & & \\
 & & 0 & 1 & \\
 & & -1 & 0 & \\
 & & & & \ddots  \end{array} \right) $$
Let
\begin{gather*}
K:=U(n,\Bbb H)=\{ g\in GL_{2n}(\C) \,|\, g^TJg=J,\ g^{-1}=\overline{g}^T\},\\
s_\theta := \mathrm{diag}(e^{i\theta},e^{-i\theta}, 1,
\ldots, 1) \in K, \quad \theta\in [0,2\pi],\\
S_\theta := \{gs_\theta g^{-1}\,|\, g\in K\} \subset K.
\end{gather*}
(here $\Bbb H$ is the division algebra of quaternions
\footnote{Quaternions appear because $U(n,\Bbb H)$ is the maximal compact 
subgroup of $Sp(2n,\Bbb C)$}.
For any $s\in S_\theta$, we denote by $p_s$ the $\omega$-orthogonal
projection of $V$ to $\mathtt{Im}(1-s)$.
                                                                                
Notice that if $s\in S_\theta$, $q\in V$, and $sq=e^{i\theta}q$,
then $sJ\oq=J\overline{sq}=e^{-i\theta}J\oq$.
Take any $q\in V$ with $q^T\oq=1$. For any $y\in V$, we have
$y = \omega(J\oq,y)q + \omega(y,q)J\oq + y'$, where $y'$ is
$\omega$-orthogonal to $q$ and $J\oq$. Let
$$ s_q(\theta) := e^{i\theta p_q},\mbox{ where }
p_q := q\oq^T + J\oq q^TJ \in i\mathtt{Lie}K.$$
Note that if $q=(1,0,\ldots,0)^T$, then
$p_q=\mathrm{diag}(1,-1,0,\ldots,0)$ and $s_q(\theta)=s_\theta$.

Let $x,y\in  V$. Now
\begin{align*}
\omega(x, p_{s_q(\theta)}y) = &
\omega(x, \omega(J\oq,y)q + \omega(y,q)J\oq) \\
= & \omega(x,q)\omega(J\oq,y)-\omega(y,q)\omega(J\oq, x) \\
= & \oq^T(yx^T-xy^T)Jq \\
= & \omega((xy^T-yx^T)\oq, q).
\end{align*}
According to Subsection \ref{SP}, we have to compute
\begin{align*}
\kappa(x,y)= & t\omega(x,y)+ \int_0^{2\pi}  
(1-e^{i\theta})(1-e^{-i\theta})
c(\theta)\Big(\int_{S_\theta}
\omega(x, p_s y)s\,ds\Big)\,d\theta\\
= & t\omega(x,y)+ \int_{|q|=1} \omega((xy^T-yx^T)\oq, q)
\Big(\int_0^{2\pi} (1-e^{i\theta})(1-e^{-i\theta})
c(\theta)e^{i\theta p_q}\,d\theta\Big)\,dq \\
= & t\omega(x,y)+ \int_{|q|=1}\omega((xy^T-yx^T)\oq, q)f(p_q)\,dq\ \in 
\Ug,
\end{align*}
where $c(\theta)$ is a finite linear combination of
$\delta_0$ and its derivatives (of even order), and $f$ is some polynomial.

Recall that we identify $\Ug$ with $S\g$ via the symmetrization map.
Let $F'_m\in S\g$ be the function on $\g$ defined by
$$F'_m(A):= \int_{|q|=1}\omega((xy^T-yx^T)\oq, q)\Tr(Ap_q)^m\,dq,
\quad A\in\g.$$
Now
\begin{gather*}
\Tr(Ap_q)=\Tr(Aq\oq^T)+\Tr(AJ\oq q^TJ) = \oq^TAq + q^TJAJ\oq \\
= \oq^TAq + \oq^TJA^TJq= 2\oq^TAq = -2\omega(AJ\oq,q).
\end{gather*}
Thus,
$$ F'_m(A)= \int_{|q|=1} \omega(B\oq,q)\omega(C\oq,q)^m\,dq $$
where $B=xy^T-yx^T$ and $C=-2AJ$.
Define $F_m\in S(\mathfrak{gl}_{2n}(\C))$ by
$$ F_m(M) :=  \int_{|q|=1} \omega(M\oq,q)^{m+1}\,dq,
\quad M\in\mathfrak{gl}_{2n}(\C).$$
Note that $dF_m|_C(B)=(m+1)F'_m(A)$.
                                                                                
We compute $F_m$ as in the previous case, by considering
$$\int_V \omega(M\oq,q)^{m+1} e^{-\zeta \oq^{T}q} \,dq.$$
Now,
\begin{align*}
\sum_{m=-1}^\infty \frac{1}{(m+1)!}\int_V \omega(M\oq,q)^{m+1}
e^{-\zeta \oq^{T}q} \,dq = & \int_V e^{-\oq^T(\zeta-M^TJ)q}\,dq\\
= & \pi^{2n}\det(\zeta-M^TJ)^{-1}\\
= & \pi^{2n}\det(\zeta+JM)^{-1}
\end{align*}
Hence, $F_m(M)$ is equal to a constant times
the coefficient of $\t^{m+1}$ in the expansion
of $\det(1+\t JM)^{-1}$ as a power series in $\t$.
Differentiating $\det(1+\t JM)^{-1}$ and evaluate at $B$, we get
$$ -\Tr\big(\t JB(1+\t JM)^{-1}\big)\det(1+\t JM)^{-1}.$$
Substituiting $B=xy^T-yx^T$, $M=C=-2AJ$, 
and replacing $\t$ by $\t/2$, we get                                  
\begin{align*}
-\t & \Tr \big( J(xy^T-yx^T)(1-\t A^T)^{-1}\big)\det(1-\t A)^{-1}/2\\
=& \t \big(\omega(x,(1-\t A)^{-1}y)-\omega(y,(1-\t A)^{-1}x)\big)
\det(1-\t A)^{-1}/2\\
=& \t \omega(x,(1-\t^2 A^2)^{-1}y)\det(1-\t A)^{-1}.
\end{align*} \qed

\section{Representation theory of continuous Cherednik algebras.} 
In this section we will consider the basics 
of representation theory of 
continuous Cherednik algebras $\hh_\kappa(G),$
where $G\subset GL(\h)$, $V=\h\oplus \h^*$. 
For simplicity we will assume that $\h$ is a faithful
representation (so $t\in \Bbb C$). 

By a representation of $\hh_\kappa(G)$ we always mean
a continuous representation, i.e., a locally finite
representation of $G$ together with a compatible action 
of $V$ satisfying the main commutation relation. 

\subsection{Dunkl operators.} 
Let $g\in GL(\h)$ be an element such that the operator $1-g$ has
rank 1. Let $f\in\C[\h]$. 
Then the function $f(u)-f(gu)$ is divisible by $\beta(u)$, 
where $\beta\in {\rm Im}(1-g^{-1*})$ is a basis element. 
Thus for any $y\in \h$ we can define the regular function 
$$
K_f(g,y,u):=\beta(y)\frac{f(u)-f(gu)}{\beta(u)}
$$
on $(S_{GL(\h)}\setminus \lbrace{1\rbrace})\times \h\times \h$
(it is independent on the choice of $\beta$). 
Since $S_{GL(\h)}$ is normal, this function extends 
to $g=1$ (by zero). Hence $K_f$ may be regarded as a
regular function on the scheme $\Psi\times \h\times \h$
(as $\Psi=S_{GL(\h)}\cap G$). 

Now let $c$ be a $G$-invariant distribution on $\Psi$. 
For any $y\in \h$, define the Dunkl operator 
$D_y$ on $\Bbb C[\h]$ by the formula
$$
(D_yf)(u)=t\partial_yf(u)+(c,K_f(?,y,u)).
$$
(Here the last term means the evaluation 
of the distribution $c$ on the function $K_f(?,y,u)$ 
of $g$ for fixed $y,u$). 

\begin{theorem}\label{dunk} Let $\kappa(x,y)=-t(x,y)-c(x,(1-g)y)$,
and $\hh_\kappa$ be the corresponding continuous Cherednik algebra.  
The assignment $x\to x$ for $x\in \h^*$, $g\to g$ (action of $g$
on $\C[\h]$), and $y\to D_y$
defines a representation of $\hh_\kappa$ on $\C[\h]$.
\end{theorem}

\begin{proof}
Let $\hh_\kappa^+$ be the subalgebra of $\hh_\kappa$ generated by
$\h$ and $\O(G)^*$. Let $\Bbb C$ be the trivial representation of
this algebra, where $y\in \h$ acts by zero, and $g\in G$ 
by $1$. Let $M:=\hh_\kappa\otimes_{\hh_\kappa^+}\Bbb C$. 
It follows from the PBW theorem for $\hh_\kappa$ that 
$M=\Bbb C[\h]$, where $x\in \h^*$ acts by $x$ and 
$g\in G$ by $g$. It remains to compute the action of $y$. 
Let $x_1,...,x_r\in \h^*$. Then we have 
$$
y\circ x_1...x_r=[y,x_1...x_r]\cdot 1=
\sum_i x_1...x_{i-1}[y,x_i]x_{i+1}...x_r\cdot 1=
$$
$$
\sum_i
(t(y,x_i)x_1...x_{i-1}x_{i+1}...x_r\cdot 1+
(c,(y,x_i-x_i^g)x_1...x_{i-1}x_{i+1}^g...x_r^g\cdot 1)=
$$ 
$$
t\partial_y (x_1...x_r)+(c,K_{x_1...x_r}(?,y,u)).
$$
(Here as usual $(c,f)$, for $f\in \O(G)$, is the evaluation
of the distribution $c$ on the function $f$).
The theorem is proved. 
\end{proof} 

The representation defined above is called the Dunkl-Cherednik
representation of the Cherednik algebra $\hh_\kappa$. 
If $G$ is finite, it is the usual Dunkl-Cherednik representation
of the rational Cherednik algebra. 

\begin{corollary} 
We have $[D_y,D_{y'}]=0$ for any $y,y'\in \h$. 
\end{corollary}                                

\subsection{Category $\mathcal O$.}       
Consider a continuous Cherednik algebra $\hh_\kappa$, with $t=1$
and arbitrary $c$. Let $\dim \h=d$. 
Like for usual Cherednik algebras, we can define the Euler element 
$$
\bold h=c+\sum_i x_iy_i+\frac{d}{2}, 
$$
where $y_i$ is a basis of $\h$
and $x_i$ the dual basis of $\h^*$. 
This element commutes with $\O(G)^*$. 

\begin{proposition} We have  
$[\bold h,x]=x$, $[\bold h,y]=-y$ for $x\in \h^*$, $y\in \h$. 
\end{proposition}

\begin{proof} We have
$$
[\bold h,x]=[c,x]+\sum_i x_i[y_i,x]=[c,x]+x+(1-g)x\cdot c=x.
$$
The second identity is proved similarly.
\end{proof}

\begin{definition}
Category $\O$ for $\hh_\kappa$ consists of finitely generated
$\hh_\kappa$-modules which are locally finite with respect to
$\bold h$, with real part of the spectrum of $\bold h$ bounded
below. 
\end{definition}

The properties of category $\O$ are as usual; let us summarize
some of them. The proofs are standard, see \cite{GGOR}. 

1. Any object of $\O$ is a finitely generated $\C[\h]$-module. 

2. For every irreducible finite dimensional $G$-module $Y$ 
we can define the Verma module 
$M(Y)=\hh_\kappa\otimes_{\hh_\kappa^+}Y$ 
(where $y\in \h$ act in $Y$ by zero). 
In particular, $M(\Bbb C)$ is the Dunkl-Cherednik
representation. The module $M(Y)$ has an 
irreducible quotient $L(Y)$, and every irreducible object in $\O$
is isomorphic to $L(Y)$. 

{\bf Remarks.}  1. In the case of Hermitian symmetric pairs
($\hh_\kappa=A(\bold G,\bold K)$), the modules $L(Y)$ 
(for appropriate $Y$) are the representations of ``holomorphic discrete series''. 

2. For $c=0$, one has $M(Y)=L(Y)$, and category $\O$ is semisimple. 
This means that this is so for Weil generic $c$, i.e. outside a countable 
union of algebraic hypersurfaces.
(An interesting problem is to find exactly for which $c$ 
this happens.) 

3. For any module $M\in \O$, we can define 
the character $\chi_M=\Tr_M(gt^{\bold h})$, 
which is a series in $t$ with coefficients in the representation
ring of $G$. For example, 
$$
\chi_{M(Y)}=\frac{t^{c_Y+d/2}}{\det|_{\h^*}(1-gt)},
$$
where $c_Y$ is the eigenvalue of $c$ in $Y$. 
One of the main problems of representation theory 
of $\hh_\kappa$ is to find the characters of $L(Y)$. 

To keep this paper short, we will not develop representation
theory of continuous Cherednik algebras systematically. 
Rather, to give the reader the flavor of this theory, 
we will study the question of reducibility of the Dunkl-Cherednik
representation for $G=GL(\h)$ and $G=O(\h)$.
\footnote{Note that for $G=GL_d$, the representation theory of
$\hh_\kappa$ is in some sense a ``deformation'' of representation
theory of the Lie group $SU(d,1)$ (with infinitely many parameters)}
 
\subsection{Representations of $\hh_\kappa$ for $G=GL(\h)$.}

\begin{proposition} 
(i) The Dunkl-Cherednik representation $M(\Bbb C)$ has singular
vectors in degree $N$ if and only if 
\begin{equation}\label{c=N}
c_{S^N\h^*}-c_{\Bbb C}=N.
\end{equation} 
Thus $M(\Bbb C)$ is irreducible iff equality (\ref{c=N}) fails for all 
$N\ge 1$. 

(ii) If $N$ is the smallest positive solution of (\ref{c=N})
then $L(\Bbb C)=\oplus_{j=0}^{N-1}S^j\h^*$ 
is the quotient of $\Bbb C[\h]$ by 
the ideal generated by polynomials of degree $N$. 
\end{proposition}

\begin{proof}
(i) The action of $y\in \h$ defines a homomorphism of $G$-modules
$\h\otimes S^N\h^*\to S^{N-1}\h^*$. There is only one such
homomorphism up to scaling, so its vanishing is equivalent to one
linear equation in $c$. If this equation is satisfied then there is a
nontrivial homomorphism $M(S^N\h^*)\to M(\C)$ of degree $N$,
which implies (\ref{c=N}). Thus the linear equation in question
is exactly (\ref{c=N}), as desired. 

(ii) is clear from the proof of (i). 
\end{proof}

\subsection{Representations of $\hh_\kappa$ for $G=O(\h)$.}

Let $G=O(\h)$, $\dim(\h)=d$. We fix an isomorphism $\h\to \h^*$
defined by the inner product on $\h$. Let $s\in O(\h)$ 
be an orthogonal reflection, and ${\rm Ad}(G)s=S\setminus \lbrace{1\rbrace}$ 
be the conjugacy class of reflections. 

In this case, it is sufficient to consider 
$c=k\Delta$, where $\Delta$ is the integral over the conjugacy
class of reflections (=invariant projection to the trivial
$G$-subrepresentation $\O(Ad(G)s)\to \Bbb C$), and $k\in \Bbb C$. Thus we will denote the
algebra $\hh_\kappa$ with $t=1$ by $\hh_k$. 

The algebra $\hh_k$ has an $sl(2)$ subalgebra, 
given by $F=\frac{1}{2}\sum x_i^2$, $E=-\frac{1}{2}
\sum y_i^2$, $H=-\frac{1}{2}(x_iy_i+y_ix_i)$
(it is easy to show that $H=-\bold h$). 

The highest eigenvalue of $H$ on the Verma module 
$M(Y)$ is 
$$
\lambda_Y(k)=-\frac{d}{2}-k\frac{\Tr|_Y(s)}{\dim(Y)}.
$$ 
This shows that in order to have a nontrivial map $M(Y)\to M(Y')$,
the number $k$ must be rational. Thus, if $k$ is not rational,
the category $\O=\O_k$ is semisimple, as $M(Y)=L(Y)$ for all $Y$. 

{\bf Definition.} The Shapovalov form $M(Y)\times M(Y^*)\to \Bbb C$ is the bilinear
$G$-invariant form coinciding with the standard one on $Y\times Y^*$, and 
satisfying the equations $(yv,w)=(v,y^*w)$, $(v,yw)=(y^*v,w)$, where $y\in \h$, and 
$y^*\in \h^*$ is the element corresponding to $y$ under the
identification $\h\to \h^*$. 

It is easy to show that the Shapovalov form exists and is
unique. Let $K(Y)$ be the kernel of the Shapovalov form, then 
$L(Y)=M(Y)/K(Y)$. 

Let us study the kernel $K(Y)$. Any representation 
of $\hh_k$ is in particular a representation of the dual pair
$(G,sl(2))$. Moreover, for generic $k$, $M(Y)=\oplus_{n\ge 0}Q(n,Y)\otimes
M_{\lambda_Y(k)-n}$, where $M_\mu$ is the highest weight Verma
module over $sl(2)$ of highest weight $\mu$, and 
$Q(n,Y)$ is the finite dimensional space of highest weight 
vectors of degree $n$ under $sl(2)$. From character
considerations, it is easy to see that  
the representation $Q(n,Y)$ of $G$ is isomorphic to 
${\rm Harm}(n)\otimes Y$, where ${\rm Harm}(n)$ is the space of harmonic 
polynomials of degree $n$. 

Let us now study the structure of the Dunkl-Cherednik representation 
$M(\C)$; i.e., $Y=\C$, $\lambda_Y(k)=:\lambda(k)=-d/2-k$. 
Let $L_{m}$ denote the finite dimensional irreducible
representation of $sl(2)$ with highest weight $m$ ($m\in \Bbb Z_+$). 

\begin{proposition}
If $k=-d/2-m$, $m\in \Bbb Z_+$, then
$L(\C)=\oplus_{n=0}^{m}L_{m-n}\otimes Harm(n)$.
Otherwise $L(\C)=M(\C)$.
\end{proposition}

{\bf Remark.} Let $d=2$, and $k=-1-m$. Then the representation
$L(\C)$ has dimension $k^2$ and is a limit of representations 
of dimension $k^2$ for finite dihedral groups $G$ 
found in Theorem 3.3.1 of \cite{Ch} (first row of the table). 

\begin{proof} 
Since all representations
${\rm Harm}(n)$ are irreducible and non-isomorphic, for any $k$ we
have $M(\C)=\oplus_{n\ge 0}M(n)\otimes {\rm Harm}(n)$, where 
$M(n)$ is an $sl(2)$-module having the same character 
as $M_{\lambda(k)-n}$. But since $F=r^2:=\sum x_i^2$, 
the module $M(\C)$ has no singular vectors under $F$. 
Hence $M(n)=M_{\lambda(k)-n}$, and thus 
$M(\C)=\oplus_{n\ge 0}M_{\lambda(k)-n}\otimes {\rm Harm}(n)$.

Thus the Shapovalov form on $M(\C)$ is the sum of the tensor
product forms on \linebreak $M_{\lambda(k)-n}\otimes Harm(n)$ 
with coefficients. 
For nonnegative integer $\lambda(k)-n$, the Shapovalov form must vanish 
on $M_{-\lambda(k)+n-2}\subset M_{\lambda(k)-n}$. 
Furthermore, we claim that all
$M_{\lambda(k)-n}\otimes {\rm Harm}(n)$ with $n>\lambda(k)$ belong to
the kernel $K(Y)$. Indeed, $K(Y)$ is a $G$-invariant proper
graded ideal in $\C[\h]$. So its zero set is either the 
cone $r^2=0$ or the origin. In particular, $r^{2n}|_{L(Y)}=0$ for
some $n$. This proves the claim. 

Thus, the zero set of $K(Y)$ is the origin, and $L(Y)$ is identified
with a $SL(2)\times O(d)$-subrepresentation of the finite dimensional 
representation $\oplus_{n=0}^{\lambda(k)}L_{\lambda(k)-n}\otimes {\rm
Harm}(n)$. 

It remains to show that this subrepresentation is in fact the entire 
representation. To do so, note that the component $L_{\lambda(k)}$ 
corresponding to $n=0$ must occur, as it contains the highest weight vector.
This implies that $r^{2p}\ne 0$ in $L(\Bbb C)$ for $p=0,...,\lambda(k)$. 
Now consider the component corresponding to some $n\le \lambda(k)$.
Consider its $sl(2)$-highest weight part, 
$Harm(n)$. This is a representation of $O(d)$ 
defined over the real numbers, so it has a real basis $f_i$ orthonormal 
under the invariant inner product. Consider the function 
$F=\sum_i f_i^2$. This function is rotationally invariant 
and generically positive on the real subspace 
(as it is a sum of squares), so it is a nonzero 
multiple of $r^{2n}$. This implies that the vector 
$\sum f_i\cdot f_i$, where the first factor $f_i$ belongs to $\hh_k$
and the second one to $L(\Bbb C)$, is nonzero. So at least one 
$f_i\in L(\Bbb C)$ is nonzero, and the entire $n$-th summand must occur 
in $L(\Bbb C)$. 

The proposition is proved. 
\end{proof}

{\bf Remark.} 
In fact, the proof of the proposition shows that 
any representation of $\hh_k$ from  the category $\O$
which does not have full support in $\h$ is finite dimensional
(as it is $sl(2)$-finite). So, the quotient of category 
of $\O$ by the category of finite dimensional representations
maps injectively into ${\rm Rep} O(n-1)$ through the fiber functor
at a generic point of $\h$. This functor is somewhat analogous to the 
KZ functor from \cite{GGOR}. 

\section{Case of   wreath-products}

\subsection{Wreath products.} Many  interesting  examples of  continuous
symplectic reflection algebras may be obtained by means of
wreath-product type construction. Specifically, let 
$\Gamma$ be a reductive subgroup 
of $SL_2$, and $G=S_n\ltimes \Gamma^n\subset Sp_{2n}$. 
When $\Gamma$ is finite (i.e., corresponds
to a simply-laced affine Dynkin diagram), the 
corresponding symplectic reflection algebras 
were considered in \cite{EG}. Thus here we will consider only 
the cases when $\Gamma$ is infinite. 
There are three such subgroups $\Gamma$: 

\begin{example} \label{eg4}
$\Gamma=\Bbb C^*=GL_1$ (the maximal torus). 
In this case $S/G=\Bbb C^*\cup \lbrace{s\rbrace}$, 
where $s$ the conjugacy class
of $s_{ij}$.\footnote{In all three examples, the class $s$ is
absent if $n=1$.}  
\end{example}

\begin{example} \label{eg5}
$\Gamma=\widetilde{O}_2$ (the normalizer of the maximal
torus). In this case, $S/G=\Bbb C^*\cup \lbrace{s,\sigma\rbrace}$, 
where $\sigma$ is a non-diagonal element of $\widetilde{O}_2$.
\end{example}

\begin{example} \label{eg6}
$\Gamma=SL_2$. 
In this case, 
$S/G=\Bbb C^*\cup \lbrace{s\rbrace}$.
\end{example}

These correspond to the infinite affine Dynkin diagrams
$A_\infty, D_\infty, A_{\infty/2}$ depicted in Figure 1
by the standard rule of McKay's correspondence (see below).  

\begin{proposition}\label{S/G}
Let $G$ be as in example \ref{eg4},\ref{eg5}, or \ref{eg6}. 
Then the most general continuous symplectic reflection algebra 
for $G$ is obtained by taking $c\in \cc(S)$.  
\end{proposition}

\begin{proof}
In Example \ref{eg4}, 
it is easy to show that $\Psi$ is reduced (i.e., $\Psi=S$), and
thus the statement follows from Corollary \ref{sra=cher}. 

In examples \ref{eg5},\ref{eg6}, 
we have $\Sigma=\Sigma_0\cup Ad(G)s, S=S_0\cup Ad(G)s$, 
where $\Sigma_0$ and $S_0$ are the intersections of $\Sigma$ and
$S$ with $G_0:=\Gamma^n$. Also, one has $\Bbb C[G_0]^{G_0}=\Bbb C[T_1,...,T_n]$, where
the functions $T_i$ are defined by the formula 
$T_i(\gamma_1,...,\gamma_n)=\Tr(\gamma_i-1)$. 
The defining equation of $\Sigma$ introduces the relations
$T_iT_j=0$, $i\ne j$. This implies that 
$\Sigma/G=S/G$, and the statement follows. 
\end{proof}

We see that in all three examples, for $n>1$ one has 
$\cc(S)=\cc(S_0)\oplus \Bbb C\Delta$, where 
$\Delta$ is the integral over $Ad(G)s$. 
Moreover $\cc(S_0)$ is naturally isomorphic to 
$\cc(\Gamma)$ (as $\Gamma$ is $S$ for $n=1$); we will identify these spaces using this
isomorphism. 

\subsection{Representations of wreath product algebras.}

The theory of wreath product algebras for infinite groups
$\Gamma$ is rather similar to the theory in the finite case. 
As an example let us work out finite dimensional 
representatations of $\hh_\kappa$ which are irreducible 
as $G$-modules. For finite groups $\Gamma$, this was done 
in \cite{M}. Since the proofs in the infinite case 
are the same as in the finite one, we will omit them. 

As before, we will denote by $\Gamma$ one of the groups
$GL_1$, $\widetilde{O}_2$ or $SL_2$, and let
$G:= S_n\ltimes \Gamma^n$. We shall write $L$ for the
2-dimensional tautological representation of $\Gamma$;
$V=L^{\oplus n}$.

Let $\{Y_i\}_{i\in I}$ be the set of finite dimensional irreducible
representations (up to isomorphism) of $\Gamma$. 
We associate a graph to $\Gamma$ in the following way:
the set of vertices of the graph is $I$, and the number of edges
joining $i,j\in I$ is the multiplicity of $Y_i$ in $L\ot Y_j$.
See Figure 1, where
$A_{\infty}$, $D_\infty$, and $A_{\infty/2}$ are the graphs
associated, respectively, to $GL_1$, $\widetilde{O}_2$,
and $SL_2$.

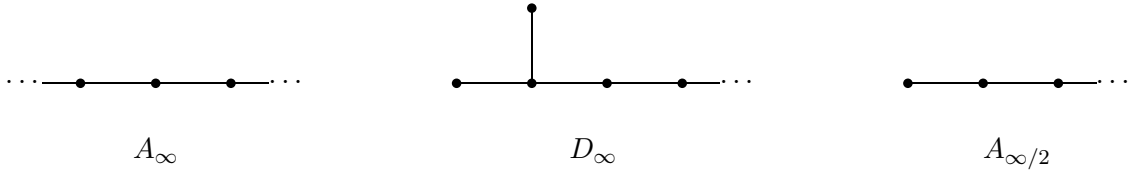
\begin{figure} 

\begin{pspicture}(0,0)(15,2 )

\dotnode(1,1){A} \dotnode(2,1){B} \dotnode(3,1){C}
\put(0,1){$\ldots$} \put(0.5,1){\line(1,0){3}}
\put(3.5,1){$\ldots$}
\put(1.7,0){$A_{\infty}$}

\dotnode(6,1){D} \dotnode(7,1){E} \dotnode(8,1){F} \dotnode(9,1){G}
\dotnode(7,2){H}
\put(6,1){\line(1,0){3.5}} 
\put(9.5,1){$\ldots$}
\put(7,1){\line(0,1){1}}
\put(7.5,0){$D_\infty$}

\dotnode(12,1){I} \dotnode(13,1){J} \dotnode(14,1){K}
\put(12,1){\line(1,0){2.5}} \put(14.5,1){$\ldots$}
\put(13,0){$A_{\infty/2}$}

\end{pspicture}
\caption{Graphs associated to $GL_1$, $\widetilde{O}_2$,
and $SL_2$.}
\end{figure}

In the theorem below, we shall use the following notations.
Let $\underline n=(n_i)_{i\in I}$, where $n_i\in\mathbb{Z}_{\geq 0}$,
and $\sum_{i\in I} n_i =n$. Let $I(\underline n):= \{i\in I\mid n_i>0\}$.
Let $Y:=\bigotimes_{i\in I(\underline n)}
Y_i^{\otimes n_i}$.
For each $i\in I(\underline n)$, we shall regard $S_{n_i}$ as the subgroup
of $S_n$ which permutes the factors in $Y_i^{\otimes n_i}$ in $Y$.
Let $S_{\underline n} := \prod_{i\in I(\underline n)} S_{n_i}
\subset S_n$.
Let $W_i$ be an irreducible representation of $S_{n_i}$,
and let $W:=\bigotimes_{i\in I(\underline n)} W_i$.
Thus, the group $S_{\underline n}\ltimes \Gamma^n$ acts on
$W\otimes Y$ (where $\Gamma^n$ acts trivially on $W$).
Denote by $W\otimes Y\uparrow$ the (irreducible) representation 
of $G$ induced from the representation $W\otimes Y$
of its subgroup $S_{\underline n}\ltimes \Gamma^n$.

We also let
$$
\kappa(x,y)=\omega(x,y)\delta_1+\omega((1-g)x,(1-g)y)(c+k\Delta),
$$ 
where $c\in \cc(\Gamma)$. We know that this gives rise to the most
general continuous symplectic reflection algebra (with $t=1$).

\begin{theorem} \label{montarani}
(i) Assume $k\neq 0$.
The representation $W\otimes Y\uparrow$ of $G$ extends to
a representation of $\hh_\kappa(G)$ where elements of
$V$ act by $0$, if and only if the following conditions
are satisfied:

{\rm (1)} For all $i\in I(\underline n)$,
the irreducible representation $W_i$ of $S_{n_i}$ has 
rectangular Young diagram, of size $a_i\times b_i$.

{\rm (2)} For all $i,j\in I(\underline n)$ (where $i\neq j$), 
the vertices $i$ and $j$ 
in the graph associated to $\Gamma$ are not adjacent.

{\rm (3)} For each $i\in I(\underline n)$, one has 
$$
\dim Y_i+2k(b_i-a_i)+(c,\chi_{Y_i}(g)\det(1-g|_L))=0,
$$
where $\chi_{Y_i}$ is the character of $Y_i$. 

(ii) If $E$ is a finite dimensional representation
of $\hh_\kappa(G)$ which is irreducible as a $G$-module, then 
$V$ acts by zero in $E$; in particular, $E$ is as 
specified in (i). 
\end{theorem}

The proof of Theorem \ref{montarani} is similar to the 
proof of \cite[Theorem 3.1, Theorem 4.1]{M}.

{\footnotesize{

\noindent
Department of Mathematics,  Massachusetts Institute 
of Technology, Cambridge, MA 02139, USA\\
{\tt etingof@math.mit.edu}
\smallskip

\noindent
Department of Mathematics,  Massachusetts Institute 
of Technology, Cambridge, MA 02139, USA\\
{\tt wlgan@math.mit.edu}
\smallskip

\noindent
Department of Mathematics, University of Chicago,
Chicago, IL 60637, USA\\
{\tt ginzburg@math.uchicago.edu}
}}

\end{document}